\documentclass[twocolumn,amssymb,aps,secnumarabic,%
    nofootinbib,superscriptaddress]{revtex4}
\usepackage{amsthm,amsmath,times,mathptmx,pstricks,pst-node,subfigure}
\newtheorem{theorem}{Theorem}[section]

\newtheorem{lemma}[theorem]{Lemma}
\newtheorem{conjecture}[theorem]{Conjecture}
\newtheorem{corollary}[theorem]{Corollary}
\newtheorem{question}[theorem]{Question}
\theoremstyle{remark}
\newtheorem*{remark}{Remark}
\newtheorem*{example}{Example}

\renewcommand{\bar}{\overline}
\renewcommand{\tilde}{\widetilde}
\renewcommand{\tensor}{\otimes}
\newcommand{\R}{{\mathbb{R}}}
\newcommand{\C}{{\mathbb{C}}}
\newcommand{\F}{{\mathbb{F}}}
\newcommand{\Q}{{\mathbb{C}}}
\renewcommand{\H}{{\mathbb{H}}}
\newcommand{\cA}{{\mathcal{A}}}
\newcommand{\cB}{{\mathcal{B}}}
\newcommand{\cC}{{\mathcal{C}}}
\newcommand{\cD}{{\mathcal{D}}}
\newcommand{\tcC}{\tilde{\mathcal{C}}}
\newcommand{\End}{\mathrm{End}}

\newcommand{\Hom}{\mathrm{Hom}}
\newcommand{\Gal}{\mathrm{Gal}}
\newcommand{\Inv}{\mathrm{Inv}}
\newcommand{\Mat}{\mathrm{Mat}}
\newcommand{\cmod}[1]{\mbox{-mod${}_{#1}$-}}
\DeclareMathOperator{\im}{im}

\newcommand{\ie}{\textit{i.e.}}

\newcommand{\thm}[1]{Theorem~\ref{#1}}
\newcommand{\prop}[1]{Proposition~\ref{#1}}
\renewcommand{\sec}[1]{Section~\ref{#1}}
\newcommand{\lem}[1]{Lemma~\ref{#1}}
\newcommand{\eq}[2]{\begin{equation}\label{#1}#2\end{equation}}

\newcommand{\que}[1]{Question~\ref{#1}}

\psset{linewidth=.35pt,dash=2.5pt 2.5pt,arrowsize=3pt 3,nodesep=.15}
\SpecialCoor

\begin{document}
\title{Finite, connected, semisimple, rigid tensor categories are linear}
\author{Greg Kuperberg}
\email[Email: ]{greg@math.ucdavis.edu}
\thanks{Supported by NSF grant DMS \#0072342}
\affiliation{Department of Mathematics,
    University of California, Davis, CA 95616}

\begin{abstract}
Fusion categories are fundamental objects in quantum algebra, but their
definition is narrow in some respects.  By definition a fusion category must be
$k$-linear for some field $k$, and every simple object $V$ is strongly simple,
meaning that $\End(V) = k$.  We prove that linearity follows automatically from
semisimplicity: Every connected, finite, semisimple, rigid, monoidal category
$\C$ is $k$-linear and finite-dimensional for some field $k$.  Barring
inseparable extensions, such a category becomes a multifusion category after
passing to an algebraic extension of $k$.

The proof depends on a result in Galois theory of independent interest, namely
a finiteness theorem for abstract composita.
\end{abstract}
\maketitle

\section{Introduction}
\label{s:intro}

We take as prerequisites to this article the first two chapters of a survey of
Bakalov and Kirillov \cite{BK:lectures} and the introductions to articles by
Etingof, Nikshych, and Ostrik \cite{ENO:fusion} and M\"uger \cite{Mueger:top1}. 
Following these three works, a \emph{fusion category} is a $k$-linear, finite,
strongly semisimple rigid tensor category. (Precise definitions of these terms
are given below.)  The previous works also present the structure
theory and applications of fusion categories.  But although fusion categories
are an important a fairly general class of objects, their definition is
narrow in some respects.

If $\cC$ is any abelian, $k$-linear category for some field $k$, we say that an
object $V \in \cC$ is \emph{strongly simple} if $\End(V) = k$. For example, if
$k$ is algebraically closed, then Schur's Lemma says that every simple object
is strongly simple.  If $\cC$ is semisimple and every simple object is strongly
simple, then we say that $\cC$ is \emph{strongly semisimple}.

A fusion category $\cC$ is endowed with a field $k$ over which it must be
linear and strongly semisimple.  It is also assumed that the identity object
$I$ is simple.  But as noted previously \cite{ENO:fusion}, it is reasonable to
drop the condition that $I$ is simple. In this case $\cC$ is a
\emph{multifusion category}, or a folded form of a fusion 2-category in which
each identity 1-morphism is simple.  (To ``fold'' an abelian
2-category means to combine finitely many objects into one
object whose 1-identity is non-simple; see \sec{s:semi}.)

The aim of this article is to show that if $\cC$ is semisimple and suitably
finite, then suitably finite linearity appears automatically.  We adopt the
natural generalization to 2-categories in the statement of the main result.

\begin{theorem} Let $\cC$ be a connected, semisimple, rigid 2-category with
finitely many types of simple 1-morphisms.  Then there exists a field $k$ over
which it is linear and all $\Hom$ spaces between 1-morphisms are
finite-dimensional.  Taking $k_s$ to be the separable closure of $k$, $k_s
\tensor \cC$ is a folded form of a semisimple 2-category over $k_s$ with
strongly simple 1-identities. If $\cC$ has no inseparable extensions of $k$,
then $k_s \tensor \cC$ is strongly semisimple.
\label{th:main} \end{theorem}

To conclude this introduction we offer two related examples,
one in which \thm{th:main} applies, and one in which it does not.

First, consider the category $\cA = \C\cmod{\R}\C$ of finite-dimensional real
vector spaces with the extra structure of bimodules over the complex numbers. 
The category $\cA$ is semisimple, and it is monoidal with respect to tensoring
in the middle.  It has two simple objects $I$ and $A$, both of which have real
dimension 2.  But in $I$, left and right complex multiplication agree, while
in $A$, they are conjugate.  The reader can check that $A \tensor A \cong I$,
from which it follows that $\cA$ is rigid, and that $\End(I) = \End(A) = \C$. 
Nonetheless, $\cA$ is not $\C$-linear.  It is $\R$-linear, while $\C \tensor_\R
\cA$ is a multifusion category which unfolds to a 2-category $\tilde{\C
\tensor_\R \cA}$ with two objects.  Each endocategory of $\tilde{\C
\tensor_\R \cA}$ is the category of complex representations of $\Gal(\C/\R)$.

Second, let $B$ be the field $\C(y)$ with a $\C(x)$-bimodule structure defined
as follows:  Left multiplication by $x$ is defined as multiplication by $y$,
while right multiplication is defined as multiplication by $y^2$.  Let $B^*$ be
$B$ with left and right switched, and let $\cB$ be the abelian monoidal
category of bimodules over $\C(x)$ generated by $B$ and $B^*$.  It is not  hard
to show that $\cB$ is semisimple and that the simple objects are $I = \C(x)$,
$B^{\tensor n}$, and $(B^*)^{\tensor n}$. In this case, the largest field over
which $\cB$ is linear is $\C$, which is already algebraically closed.  No
change of base field of $\cB$ renders it strongly semisimple.

\acknowledgments

The author would like to thank Michael Khovanov and especially Bjorn
Poonen for helpful discussions.

\section{Semisimple categories}
\label{s:semi}

We assume various relevant definitions from Mac Lane \cite{MacLane:gtm}
and M\"uger \cite{Mueger:top1}:
additive, abelian, monoidal, $k$-linear, etc.  An object in an abelian category is
\emph{simple} if it has no subobjects.  An abelian category is
\emph{semisimple} if every object is a direct sum of finitely many simple
objects.  As mentioned in the introduction, if $V$ is a simple object in an
abelian, $k$-linear category, then it is \emph{strongly simple} if $\End(V) =
k$.  A semisimple, $k$-linear category is \emph{strongly semisimple} if all
simple objects are strongly simple.

A (strict) monoidal category $\cC$ can be reinterpreted as a 2-category $\cC'$
with one object.  This phenomenon is known as \emph{dimension shifting},
because $n$-morphisms in $\cC$ become $(n+1)$-morphisms in $\cC'$. In light of
this relationship, we will use $V \tensor W$ to denote the composition of
1-morphisms $V \in \Hom_1(A,B)$ and $W \in \Hom_1(B,C)$, and $f \tensor g$ for
the attendant ``horizontal'' composition of 2-morphisms.  We use $fg$ or $f
\circ g$ to denote ``vertical'' composition of two 2-morphisms $f$ and $g$ in
the same category $\Hom_1(A,B)$.  A 2-category is also additive, abelian,
$k$-linear, semisimple, or strongly semisimple if each category $\Hom_1(A,B)$
has the same property.  If additivity or linear is part of the structure, we
also assume that $f \tensor g$ is biadditive or bilinear in $f$ and $g$.

\begin{remark}
It is common to assume, at least intuitively, that categories are
\emph{skeletal} (there is only one object of each isomorphism type), because
every category can be made skeletal. It is also common to assume that monoidal
categories are \emph{strict} ($\tensor$ is strictly associative), because every
monoidal category can be made strict.  But most monoidal categories cannot be
made simultaneously strict and skeletal!  This is the origin of the non-trivial
structure of associators. On balance, we prefer strictness and we will not
assume that categories are skeletal.  However, we can assume that 2-categories
are skeletal at the level of objects.
\end{remark}

Let $\cC$ be a semisimple monoidal category and suppose that the identity
object $I$ is not simple.  Then
$$I = \bigoplus_{A \in S} I_A,$$
where $S$ is some indexing set and each $I_A$ is simple. As previously
noted \cite{ENO:fusion}, $I$ is necessarily multiplicity-free, \ie, $I_A
\not\cong I_B$ when $A \ne B$, and $\cC$ can be reorganized as a 2-category
$\tcC$ whose objects are the elements of $S$. In $\tcC$, the identity of $A$
is $I_A$, and the Hom space $\Hom_1(A,B)$ consists of those objects $V$ of
$\cC$ such that
$$I_A \tensor V \tensor I_B = V \in \cC.$$
(Note that $I_A \tensor V \tensor I_B$ is always a subobject of $V$.)  Thus
$\tcC$ has simple 1-identities. We call it the \emph{unfolded} form of $\cC$.

The same construction applies if $\cC$ is a semisimple 2-category such that not
all 1-identities are simple. If each 1-identity $I_A$ decomposes as
$$I_A = \bigoplus_{B \in S_A} I_B,$$
then $S = \bigcup_A S_A$ is the set of objects of the unfolded category $\tcC$.

\begin{lemma} If $V \in \Hom_1(A,B)$ is a simple 1-morphism in an abelian
2-category $\cC$, then $\End_2(V)$ is a division ring.  
\label{l:schur} \end{lemma}

\begin{proof} The lemma is a form of Schur's Lemma. Suppose
that $V$ is simple and that $f \in \End_2(V)$. Then both $\ker f$ and $\im f$
are either $0$ or $V$. If $\ker f = V$ or if $\im f = 0$, then $f = 0$.
Otherwise, if $\ker f = 0$ and $\im f = V$, then $f$ has an inverse on each
side and therefore a two-sided inverse.  Thus every non-zero $f$ has a
reciprocal and $\End_2(V)$ is a division ring.
\end{proof}

\begin{lemma} The division ring $\End_2(I_A)$ is a field.  If $V \in
\Hom_1(A,B)$ is any 1-morphism, then $\End_2(I_A)$ and $\End_2(I_B)$ embed in
the center $Z(\End_2(V))$.
\label{l:center} \end{lemma}
\begin{proof} The identity $V = I_A \tensor V$ induces a unital ring
homomorphism $\End_2(I_A) \to \End_2(V)$, which must be an inclusion since the
domain is a division ring. To show that $\End_2(I_A)$ lies in the center of
$\End_2(V)$, let $f \in \End_2(I_A)$ and $g \in \End_2(V)$.  Then
$$fg = (f \tensor 1_V)(1_{I_A} \tensor g)
    = (1_{I_A} \tensor g)(f \tensor 1_V) = gf.$$
In particular, if $V = I_A$, then $\End_2(I_A)$ lies in the center of itself,
so it is a field.
\end{proof}

\begin{remark} If $A$ is an object in any 2-category $\cC$, then
$\End_2(I_A)$ is commutative.  This fact is familiar in topology as the
commutativity of the second homotopy group $\pi_2(X)$ of a topological space
$X$.  Our proof is the usual one.
\end{remark}

Note that $\End_2(V)$ may not be commutative even if $V$ is
simple.  For example, if $V$ is the defining representation of
$\mathfrak{sl}(2,\C)$ viewed as a 4-dimensional real representation, then
$\End_2(V) = \H$, the quaternions.

In light of \lem{l:center}, let $k_A = \End_2(I_A)$, let $k_V$ be
the compositum of the inclusions $k_A$ and $k_B$ in $Z(\End(V))$, and denote the
restrictions of these inclusions to $k_V$ by
$$k_A \stackrel{\alpha_{A,V}}{\longrightarrow}
    k_V \stackrel{\alpha_{V,B}}{\longleftarrow} k_B$$
Thus $k_V$ is an \emph{abstract compositum} of the fields
$k_A$ and $k_B$.

Now let $\cC$ be a semisimple 2-category with simple 1-identities.   Then the
fields $k_A$ are entirely unrelated on different connected components of
$\cC$.  Even when $\cC$ is connected, the fields $k_A$ may differ for
different $A \in \cC$, although they must have the same 
characteristic since they are connected by abstract composita.

If $A,B \in \cC$ are objects and $V,W \in \Hom_1(A,B)$ are two 1-morphisms
connecting them, then $\Hom_2(V,W)$ is a bimodule over the two fields $k_A$ and
$k_B$.  If $k_A \ncong k_B$, then the left and right module structures
certainly differ.  But even if $\cC$ only has one object $A$, so that it is a
dimension-shifted monoidal category, the left and right $k_A$-module structures
may differ; an example was given in \sec{s:intro}.

\begin{example}  Let $\cA$ be a 2-category with two objects, $\R$ and $\C$,
and define its Hom categories by
$$\Hom_1(k_1,k_2) = k_1\cmod{\R} k_2$$
for every $k_1,k_2 \in \{\R,\C\}$.  Then $\End_2(I_k) = k$, so the
endomorphism fields of the 1-identities differ.
\end{example}

\section{Rigidity}
\label{s:rigid}

A monoidal category is \emph{rigid} if every object $V$ has both a left
dual ${}^*V$ and a right dual $V^*$ together with morphisms
\begin{align*}
a_V:I &\to V \tensor V^* & b_V:V^* &\tensor V \to I \\
c_V:I &\to {}^*V \tensor V & d_V:V &\tensor {}^*V \to I
\end{align*}
that satisfy the compatibility conditions
\begin{align*}
(1_V \tensor b_V)(a_V \tensor 1_V) &= 1_V \\
(d_V \tensor 1_V)(1_V \tensor c_V) &= 1_V.
\end{align*}
these definitions generalize readily to 2-categories.  To be explicit, if $V
\in \Hom_1(A,B)$, then  $V^*,{}^*V \in \Hom_1(B,A)$, and
\begin{align*}
a_V:I_A &\to V \tensor V^* & b_V:V^* &\tensor V \to I_B \\
c_V:I_B &\to {}^*V \tensor V & d_V:V &\tensor {}^*V \to I_A,
\end{align*}
and the compatibility conditions are the same. As explained by M\"uger
\cite{Mueger:top1}, duals can also be called \emph{adjoints}, with the
significant consequence that $V \mapsto V^*$ and $V \mapsto {}^*V$ are
contravariant endofunctors of the 2-category $\cC$.  In fact they
can be made inverse to each other, so that $V = {}^*V^*$. 

Note also that if $\cC$ is semisimple, then $V \cong V^{**}$ \cite{ENO:fusion},
although unless $\cC$ is \emph{pivotal}, these isomorphisms are not
functorial.  We will not need pivotal structure in this article.

Henceforth let $\cC$ be a semisimple, rigid 2-category with simple
1-identities.

\begin{lemma} If $V \in \Hom_1(A,B)$, then there is an isomorphism
$\sigma_V:k_V \to k_{V^*}$ that makes the following diagram commute:
$$\pspicture(-2,-1.5)(2,1.5)
\rput(-2,0){\rnode{a}{$k_A$}} \rput(0,1){\rnode{b}{$k_V$}}
\rput(0,-1){\rnode{c}{$k_{V^*}$}} \rput(2,0){\rnode{d}{$k_A$}}
\ncline{->}{a}{b} \Aput{$\alpha_{A,V}$}
\ncline{->}{a}{c} \Bput{$\alpha_{V^*,A}$}
\ncline{->}{b}{c} \Aput{$\sigma_V$}
\ncline{->}{d}{b} \Bput{$\alpha_{V,B}$}
\ncline{->}{d}{c} \Aput{$\alpha_{B,V^*}$}
\endpspicture$$
\label{l:dual} \end{lemma}
\begin{proof}  Since $V \mapsto V^*$ extends to an anti-automorphism
of $\cC$, it produces $\sigma_V$.  It is only necessary to
check that the anti-automorphism is the identity on $k_A = \End_2(I_A)$.
\end{proof}

\begin{lemma} If $V \in \Hom_1(A,B)$ is a 1-morphism in $\cC$, then $\End_2(V)$
is finite-dimensional as a left $k_A$-module and as a right $k_B$-module.
\label{l:fin}
\end{lemma}

\begin{proof} If $W \in \Hom_1(A,A)$, then it has an \emph{invariant
space} defined as
$$\Inv(W) = \Hom_2(I_A,W).$$
Since $\cC$ is semisimple, $\Inv(W)$ is finite-dimensional as a vector space
over $k_A$. If $W = V \tensor {}^*V$, then
$$\Inv(V \tensor {}^*V) \cong k_A^n$$
with $n>0$ by the existence of $d_V$. Let $\kappa_V$ be the composition of the
maps
$$\End_2(V) \to \End_2(V \tensor {}^*V) \to \End(\Inv(V \tensor {}^*V))
    \cong \Mat_n(k_A),$$
where the first term is given by  $f \mapsto f \tensor 1_{{}^*V}$ and the
second by isotypic decomposition of $V \tensor {}^*V$. The map $\kappa_V$ is
both a unital ring homomorphism and morphism of $k_A$-linear spaces.  Since
the domain of $\alpha$ is a division ring and $1 \ne 0$ in the target,
$\kappa_V$ is injective.  Since the target is a finite-dimensional vector
space over $k_A$, the domain $\End_2(V)$ is finite-dimensional as well.

The same proof works on the other side, replacing $\kappa_V$ with
$$\End_2(V) \to \End_2(V^* \tensor V) \to \End(\Inv(V^* \tensor V))
    \cong \Mat_n(k_A),$$
which for later use we call $\lambda_V$.
\end{proof}

\begin{remark} In fact, $\Inv(V \tensor {}^*V)$ and $\End_2(V)$ are isomorphic
as $k_A$-vector spaces.
\end{remark}

Combining Lemmas~\ref{l:schur} and \ref{l:fin}, the abstract
compositum $k_V$ associated to a simple 1-morphism $V \in
\Hom_1(A,B)$ is bifinite:  a finite extension of both $k_A$ and $k_B$.

\begin{remark} \thm{th:main} postulates a common finite-index subfield of all
$k_A$ and all $k_V$ on which every $\alpha_{A,V}$ is the identity.  Given that
there are only finitely many simple $V$ up to isomorphism, the fact that $k_V$
is a bifinite compositum of $k_A$ and $k_B$ suggests looking at their
intersection in $k_V$.  Unfortunately, if $k_1$ and $k_2$ are two finite-index
subfields of a field $k_3$, it does not follow that $k_1 \cap k_3$
is finite index in $k_3$.  For example, let
$$k_3 = \C(q) \qquad k_1 = \C(q^2) \qquad k_2 = \C((q-1)^2).$$
The reader can check that $k_1 \cap k_2 = \C$.  So it is not enough
to know that bifinite abstract composita connect every pair $k_A$ and $k_B$.
\end{remark}

\begin{lemma} Let $V \in \Hom_1(A,B)$ and $W \in \Hom_1(B,C)$. Then the map
$$\sigma_{V,W}:\End_2(V) \tensor_{k_B} \End_2(W) \to \End_2(V \tensor W)$$
is injective.
\label{l:inject} \end{lemma}
\begin{proof} Consider the commutative diagram
\vspace{.5cm}$$\begin{array}{c@{\hspace{1cm}}c}
\rnode{a}{\End_2(V) \tensor \End_2(W)} & \rnode{b}{\End_2(V \tensor W)} \\[1.5cm]
 & \rnode{c}{\End_2(V^* \tensor V \tensor W \tensor {}^*W)} \\[1cm]
\rnode{d}{\begin{array}{c}\End(\Inv(V^* \tensor V)) \\ \tensor \\ \End(\Inv(W \tensor {}^*W))
    \end{array}} & 
\rnode{e}{\End(\Inv(V^* \tensor V \tensor W \tensor {}^*W))}
\end{array}.$$
\ncline{->}{a}{b} \Aput{$\sigma_{V,W}$}
\ncline{->}{b}{c}
\ncline{->}{a}{d} \Aput{$\lambda_V \tensor \kappa_W$}
\ncline{->}{c}{e}
\ncline{->}{d}{e}
Since the maps $\kappa_W$ and $\lambda_V$
from \lem{l:fin} are both injective, the left
arrow is injective.  The bottom arrow is trivially injective.
Therefore the top arrow, $\sigma_{V,W}$, is also injective.
\end{proof}

\begin{lemma} If $V \in \Hom_1(A,B)$, then the abstract compositum $k_V$
is a separable extension of both $k_A$ and $k_B$.
\label{l:sep} \end{lemma}
\begin{proof} By abuse of notation, we omit the embeddings $\alpha_{A,V}$ and
$\alpha_{V,B}$. (This already arises in the statement of the lemma.)  Applying
\lem{l:inject} to $V \tensor {}^*V$, we know that $k_V \tensor_{k_B} k_V$
embeds in $\End_2(V \tensor {}^*V)$. Moreover, the subalgebra $R$ generated by
both copies of  $k_A$ in $k_V \tensor_{k_B} k_V$ lies in the center of
$\End_2(V \tensor {}^*V)$, since the latter is an algebra over $k_A$ on both
the left and the right.  It suffices to show that $R$ has nilpotent elements
when $k_V$ is inseparable over $k_B$, because this would violate the
semisimplicity of $V \tensor {}^*V$.

Let $p \ne 0$ be the common characteristic of $k_A$, $k_B$, and $k_V$. Suppose
as a special case that $k_V$ is a  non-trivial purely inseparable extension of
$k_B$ with exponent $e$.  Since $k_V$ is the compositum of $k_A$ and $k_B$,
$k_A$ contains an element $x$ which is not in $k_B$.  In this case 
$$x \tensor 1 - 1 \tensor x \ne 0 \in k_V \tensor_{k_B} k_V,$$
while
$$(x \tensor 1 - 1 \tensor x)^{p^e} = 0.$$
Thus $x \tensor 1 - 1 \tensor x$ is the desired nilpotent
element in $R$.

In this general case, $k_B$ has a maximal separable extension $s_V$ in $k_V$. 
Then $k_V \tensor_{k_B} k_V$ surjects onto $k_V \tensor_{s_V} k_V$.  Replacing
$k_B$ by $s_V$ in the previous paragraph, the image $R'$ of $R$ in $k_V
\tensor_{s_V} k_V$ has a nilpotent element.  At the same time, $k_V
\tensor_{k_B} k_V$ is a finite-dimensional algebra with respect to its left
$k_A$ structure;
$$\dim_{k_A} k_V \tensor_{k_B} k_V
= (\dim_{k_A} k_V)(\dim_{k_B} k_V) < \infty$$
by \lem{l:fin}. Thus $R$ is also a finite-dimensional (commutative) algebra
over $k_A$. Since its quotient $R'$ has a nilpotent element, $R$ must have a
nilpotent element as well.
\end{proof}

\begin{lemma} Let $V \in \Hom_1(A,B)$ and $W \in \Hom_1(B,C)$. Then any
compositum of $k_A$ and $k_C$ which occurs as a subring of $k_V
\tensor_{k_B} k_W$ is $k_X$ for some summand $X \subseteq V \tensor W$.
\label{l:summand}
\end{lemma}
\begin{proof} Since the extensions $k_V$ and $k_W$ are separable over
$k_B$,
$k_V \tensor_{k_B} k_W$ is semisimple and decomposes as a direct sum
of fields:
$$k_V \tensor_{k_B} k_W = \bigoplus_{i=1}^n k_i.$$
Let $P_i \in k_V \tensor_{k_B} k_W$ be the projection onto the summand $k_i$.
Since $k_V \tensor_{k_B} k_W$ embeds in $\End_2(V \tensor W)$, we can view
$P_i$ as a non-zero idempotent in $\End_2(V \tensor W)$ as well.  Then $\im P_i
\in \Hom_1(A,C)$ and $k_i \subseteq \End_2(\im P_i)$. It follows that the
subfield of $k_i$ generated by $\alpha_{A,V}(k_A)$ and $\alpha_{W,C}(k_C)$ is
also the field $k_X$ for any simple summand $X$ of $\im P_i$.
\end{proof}

\section{Galois theory}
\label{s:galois}

In this section we complete the proof of \thm{th:main} using a result in Galois
theory which is of separate interest.

\sec{s:rigid} constructs, from a suitable 2-category $\cC$, a finite collection
of fields $\{k_A\}$ and a finite collection of bifinite, biseparable
abstract composita
$$k_A \stackrel{\alpha_{A,V}}{\longrightarrow}
    k_V \stackrel{\alpha_{V,B}}{\longleftarrow} k_B$$
Moreover,
$k_{I_A} = k_A$ and $\alpha_{A,I_A}$ is the identity.  By \lem{l:dual},
$\alpha_{A,V^*} = \alpha_{V,A}$.  And by \lem{l:summand}, given two
abstract composita
$$k_A \stackrel{\alpha_{A,V}}{\longrightarrow}
    k_V \stackrel{\alpha_{V,B}}{\longleftarrow} k_B
    \stackrel{\alpha_{B,W}}{\longrightarrow}
    k_W \stackrel{\alpha_{W,C}}{\longleftarrow} k_C,$$
every compositum of $k_A$ and $k_C$ that appears in $k_V
\tensor_{k_B} k_W$ is $k_X$ for some summand $X \subseteq V \tensor W$.
We call this method of producing $k_X$ from $k_V$ and $k_W$
\emph{amalgamation} of composita.

\begin{theorem} Let $K = \{k_A\}$ be a finite set of fields, and let $E =
\{k_V\}$ be a finite set of biseparable, bifinite composita with embeddings 
$$k_A \stackrel{\alpha_{A,V}}{\longrightarrow}
    k_V \stackrel{\alpha_{V,B}}{\longleftarrow} k_B.$$
Suppose that $E$ contains identities, is closed with respect to duality and
amalgamation, and connects every pair of elements of $K$.  Then there is a
field $k$ and finite-index embeddings $\{\beta_A\}$ and $\{\beta_V\}$ that
form commutative triangles:
$$\pspicture[.45](-2.5,-1.5)(1,1.5)
\rput(-2,0){\rnode{a}{$k$}} \rput(0,1){\rnode{b}{$k_A$}}
\rput(0,-1){\rnode{c}{$k_V$}}
\ncline{->}{a}{b} \Aput{$\beta_A$} \ncline{->}{a}{c} \Bput{$\beta_V$}
\ncline{->}{b}{c} \Aput{$\alpha_{A,V}$}
\endpspicture$$
\label{th:galois} \end{theorem}

\prop{th:galois} can also be reformulated as the following corollary.

\begin{corollary} Let $K = \{k_A\}$ be a finite set of fields, let
$K\mbox{-Mod-}K$ be the 2-category of bimodules over elements of $K$, and let
$\cD$ be a full, connected, rigid, semisimple sub-2-category of
$K\mbox{-Mod-}K$ with finitely many simple 1-morphisms and without inseparable
extensions.  Then $\cD$ admits a forgetful functor to the category $k$-mod of
finite-dimensional vector spaces over a field $k$ which embeds in every $k_A$.
\label{c:galois}
\end{corollary}

By \lem{l:sep}, the semisimplicity of the 2-category $\cD$ eliminates the
possibility of inseparable extensions. Although $\cD$ satisfies the hypotheses
of \thm{th:main}, and although it is constructed from the more general
2-category $\cC$, we do not know a natural functor from $\cC$ to $\cD$.  

\begin{proof} Since the fields $k_A \in K$ are all connected by composita,
they all have isomorphic separable closures. Let $\Omega$ be
a field in this isomorphism class, and realize each $k_A$
arbitrarily as a subfield of $\Omega$. Let $G_A = \Gal(\Omega/k_A)$
be the absolute Galois group of $k_A$.  Finally let $\F$ be
the characteristic field of $\Omega$, either $\F_p$ or $\Q$.

If $k_V \in E$ is an extension of $k_A$, we can position $k_V$ so that $k_A
\subseteq k_V \subseteq \Omega$; the embedding $\alpha_{A,V}$ is then the
inclusion map.  But having chosen this position for $k_V$, we cannot require
that the other embedding $\alpha_{V,B}$ is the inclusion map.  Rather $k_V$
only contains a subfield isomorphic to $k_B$ and $\alpha_{V,B}$ is the
isomorphism. Let $\phi_V \in \Gal(\Omega/\F)$ be an extension to all of
$\Omega$ of the map $\alpha_{V,B}$, so that
\eq{e:gv}{G_V \stackrel{\mathrm{def}}{=} \Gal(\Omega/k_V) =
    G_A \cap \phi_V G_B \phi_V^{-1}.}
Note that, having fixed $k_A,k_B \subseteq \Omega$, the connecting
automorphism $\phi_V$ can be replaced by any other element of the double coset
$G_A \phi_V G_B$.  The double coset determines the mutual extension $k_V$ up to
its position in $\Omega$.

Next consider the tensor product $k_V \tensor_{k_B} k_W$, which, as in the
proof of \lem{l:summand}, is a direct sum of fields:
$$k_V \tensor_{k_B} k_W = \bigoplus_{i=1}^n k_i.$$
Any summand $k_i$
contains a copy of the field embeddings
$$k_A \stackrel{\alpha_{A,V}}{\longrightarrow}
    k_V \stackrel{\alpha_{V,B}}{\longleftarrow} k_B
    \stackrel{\alpha_{B,W}}{\longrightarrow}
    k_W \stackrel{\alpha_{W,C}}{\longleftarrow} k_C.$$
As discussed above, the relative position of $k_A$ and $k_B$ in $k_V$ is
described by an element of the double coset $G_A \phi_V G_B$.  Likewise the
relative position of $k_B$ and $k_C$ is described by an element of the double
coset $G_B \phi_W G_C$. Therefore the relative position of $k_A$ and $k_C$ in
the summand $k$ is given by an element $\phi$ of their product
$$G_A \phi_V G_B \phi_W G_C.$$
By hypothesis, the compositum $k_X$ of $k_A$ and $k_C$ in $k_i$ which is in $E$
and is represented by its own double coset $G_A \phi_X G_C$.  Thus
$$\phi \in G_A \phi_X G_C \subseteq G_A \phi_V G_B \phi_W G_C.$$
At the same time, universality of tensor products implies that if $\phi$ is any
element of $G_A \phi_V G_B \phi_W G_C$, the corresponding relative position of
$k_V$ and $k_W$ is represented by some summand $k \subseteq k_V \tensor_{k_B}
k_W$.

Thus the decomposition of $k_V \tensor_{k_B} k_W$ yields a decomposition of
double cosets
\eq{e:vw}{G_A \phi_V G_B \phi_W G_C
    = \bigcup_{X \in E_V} G_A \phi_X G_C}
for some subset $E_V \subseteq E$.  In addition, the duality hypothesis implies
that we can take we can take $\phi_{V*} = \phi_V^{-1}$ for some $V^* \in E$,
while the identity hypothesis implies that we can take $\phi_{I_A} = 1$ for
some $I_A \in E$.  Combining all of these facts, if $E_{A,B}$ is the set of all
mutual extensions of $k_A$ and $k_B$ in $E$, then the union of double cosets
$$H_A = \bigcup_{V \in E_{A,A}} G_A \phi_V G_A$$
is a group:  It is closed under multiplication and inversion and contains
the group $G_A$ (and therefore 1).

We claim that $G_A$ is a finite-index subgroup of $H_A$.  By
equation~\eqref{e:gv}, the number of right cosets of $G_A$ in the double coset
$G_A \phi_V G_A$ is the same as the index $[G_A:G_V]$ of $G_V$ in $G_A$, which
by hypothesis is finite.  Moreover, $H_A$ is a finite union of such double
cosets, since $E$ is finite.  This establishes the claim.  As the remaining
arguments indicate, the claim is the heart of the proof of \thm{th:galois}.

Let $f_A$ be the fixed field of $H_A \subseteq \Gal(\Omega/\F)$. Then $k_A$ is
a finite, separable extension of $f_A$, because its Galois group $G_A$ is a
finite-index subgroup of $H_A$. Moreover, $f_A$ and $f_B$ are canonically
isomorphic. Any $\phi_V \in \Gal(\Omega/\F)$ is an isomorphism between them. 
By equation~\eqref{e:vw}, any two choices for $\phi_V$ (allowing $V$ to vary as
well) differ by an element of $H_A$ on the left, and therefore all agree after
restriction to $f_A$.  Equation~\eqref{e:vw} also shows that the isomorphisms
between $f_A$, $f_B$, and $f_C$ form a commutative triangle.  Thus we can let
$k$ be a field isomorphic to all of them by maps
$$\beta_A:k \to f_A;$$
the same map $\beta_A$ can also be taken as an embedding of $k$ in $k_A$.  By
construction these field embeddings also extend to commutative triangles
$$\pspicture[.45](-2.5,-1.5)(1,1.5)
\rput(-2,0){\rnode{a}{$k$}} \rput(0,1){\rnode{b}{$k_A$}}
\rput(0,-1){\rnode{c}{$k_V$}}
\ncline{->}{a}{b} \Aput{$\beta_A$} \ncline{->}{a}{c} \Bput{$\beta_V$}
\ncline{->}{b}{c} \Aput{$\alpha_{A,V}$}
\endpspicture,$$
as desired.
\end{proof}

\section{Questions}
\label{s:questions}

\begin{question} If $\cC$ is $k$-linear over a separably closed field $k$ and
every 1-identity is strongly simple, can some $\End_2(V)$ be
a non-trivial inseparable extension of $k$?
\label{q:sep} \end{question}

It is noteworthy that if $f$ is an inseparable finite extension of $k$, then $f
\tensor_k f$ is not semisimple; this is a weak form of \lem{l:sep}.  If a
2-category $\cC$ did satisfy \que{q:sep}, then the result of base change
$\bar{k} \tensor_k \cC$ would not be multifusion because it would not be
semisimple.  This could be taken as a loophole in the structure theory of
fusion categories in characteristic $p$: such a category $\cC$ would be
``morally'' but not ``technically'' fusion.

If $\cC$ is a semisimple 2-category with simple 1-identities, then it may be
\emph{weakly right-rigid} in the sense that for every $V \in \Hom(A,B)$, there
exists $V^* \in \Hom(B,A)$ such that
$$\Inv(V \tensor V^*) \ne 0.$$
If $\cC$ is weakly right-rigid, then the structure of
$V \tensor V^* \tensor V^{**}$ induces a map
\begin{multline*}
s_V:\End_2(V) \tensor_{k_B} \Inv(V^* \tensor V^{**}) \\
    \to \Inv(V \tensor V^*) \tensor_{k_A} \End_2(V^{**}).
\end{multline*}
We can always take $V^*$ and $V^{**}$ to be simple.  If $V$ is also simple,
then either $s_V$ vanishes, or $V^{**} \cong V$ and $V^*$ is (up to
isomorphisms) both a left dual and a right dual of $V$.

\begin{conjecture} Every finite, weakly rigid, semisimple 2-category
with simple 1-identities is rigid.
\end{conjecture}

\begin{example} The representation category of $U_q(\mathfrak{sl}(2))$ becomes
weakly rigid but not rigid in the crystal limit $q \to 0$.  However, it has
infinitely many simple objects.
\end{example}

% \bibliography{qa,books}

\begin{thebibliography}{1}

\bibitem{BK:lectures}
Bojko Bakalov and Jr. Alexander~Kirillov, \emph{Lectures on tensor categories
  and modular functors}, American Mathematical Society, Providence, RI, 2001.

\bibitem{ENO:fusion}
Pavel Etingof, Dmitri Nikshych, and Viktor Ostrik, \emph{On fusion categories},
  \mbox{arXiv:math.QA/0203060}.

\bibitem{MacLane:gtm}
Saunders {Mac Lane}, \emph{Categories for the working mathematician}, second
  ed., Springer-Verlag, New York, 1998.

\bibitem{Mueger:top1}
Michael M\"uger, \emph{From subfactors to categories and topology {I}.
  {Frobenius} algebras in and {Morita} equivalence of tensor categories},
  \mbox{arXiv:math.CT/0111204}.

\end{thebibliography}

\providecommand{\bysame}{\leavevmode\hbox to3em{\hrulefill}\thinspace}

\end{document}